\documentclass[submission]{eptcs}

\usepackage{amsmath,amsfonts,amssymb,stmaryrd}

\usepackage{amsthm}

\theoremstyle{definition}

\newtheorem{example}{Example}

\usepackage{breakurl}

\usepackage{tikzit}
\usetikzlibrary{matrix}
\tikzstyle{circ}=[fill=white, draw=black, shape=circle]
\tikzstyle{blank}=[fill=white, draw=white, shape=circle]
\tikzstyle{none}=[fill=none, draw=none]
\tikzstyle{copy}=[fill=white, draw=black, shape=circle, minimum height=0.2cm, inner sep=0]
\tikzstyle{varCopy}=[fill=black, draw=black, shape=circle, minimum height=0.2cm, inner sep=0]
\tikzstyle{copy2}=[fill=black, draw=black, shape=circle, minimum height=0.2cm, inner sep=0]
\tikzstyle{1morph1}=[fill=white, draw=black, shape=rectangle, minimum width=1cm, minimum height=1cm]
\tikzstyle{1morph}=[fill=white, draw=black, shape=rectangle, minimum width=0.75cm, minimum height=0.75cm, inner sep=0.1cm]
\tikzstyle{2morph2}=[fill=white, draw=black, shape=rectangle, minimum width=1cm, minimum height=2cm]
\tikzstyle{2morph}=[fill=white, draw=black, shape=rectangle, minimum width=1cm, minimum height=1.25cm, inner sep=0.1cm]
\tikzstyle{nmorph}=[fill=white, draw=black, shape=rectangle, minimum height=6cm, minimum width=1cm, inner sep=0.1cm]
\tikzstyle{1state}=[fill=white, draw=black, regular polygon, regular polygon sides=3, minimum height=0.5cm, regular polygon rotate=-30]
\tikzstyle{dbox}=[fill=white, draw=black, dashed, shape=rectangle, minimum width=2cm, minimum height=1cm, inner sep=0.1cm]
\tikzstyle{vdbox}=[fill=white, draw=black, dashed, shape=rectangle, minimum width=2cm, minimum height=1.5cm, inner sep=0.1cm]
\tikzstyle{bigbox}=[fill=white, draw=black, dashed, shape=rectangle, minimum width=2cm, minimum height=4cm, inner sep=0.1cm]
\tikzstyle{2state}=[inner sep=0.05cm, fill=white, draw=black, isosceles triangle, minimum width=1.25cm, isosceles triangle apex angle=90, shape border rotate=180]
\tikzstyle{var2state}=[inner sep=0.05cm, fill=white, draw=black, isosceles triangle, minimum width=1.25cm, isosceles triangle apex angle=60, shape border rotate=180]
\tikzstyle{g2state}=[inner sep=0.05cm, fill=white, draw=black, isosceles triangle, minimum width=6cm, isosceles triangle apex angle=110, shape border rotate=180]
\tikzstyle{bigstate}=[inner sep=0.05cm, fill=white, draw=black, isosceles triangle, minimum width=3cm, isosceles triangle apex angle=110, shape border rotate=180]
\tikzstyle{bigeffect}=[inner sep=0.05cm, fill=white, draw=black, isosceles triangle, minimum width=3cm, isosceles triangle apex angle=110]
\tikzstyle{g2effect}=[inner sep=0.05cm, fill=white, draw=black, isosceles triangle, minimum width=6cm, isosceles triangle apex angle=110]
\tikzstyle{2effect}=[inner sep=0.05cm, fill=white, draw=black, isosceles triangle, minimum width=1.25cm, isosceles triangle apex angle=90]
\tikzstyle{b2effect}=[inner sep=0.05cm, fill=white, draw=black, isosceles triangle, minimum width=2cm, isosceles triangle apex angle=90]
\tikzstyle{midArrow}=[-, decoration={{markings,mark=at position .5 with {\arrow{>}}}}, postaction=decorate]

\tikzstyle{Medium box}=[fill=white, draw=black, shape=rectangle, tikzit shape=rectangle, minimum width=1.5cm, minimum height=1.5cm]
\tikzstyle{Black arrow}=[->]
\tikzstyle{Gray line}=[-, draw={rgb,255: red,191; green,191; blue,191}, line width=0.8]

\tikzstyle{arrow}=[->]

\newcommand\R{\mathbb R} % real numbers
\newcommand\RR{\mathbb R} % also real numbers
\newcommand\C{\mathcal C} % generic category
\newcommand\D{\mathcal D} % another generic category
\newcommand\M{\mathcal M} % generic residual category
\newcommand\E{\mathbb E} % expectation map

\newcommand\op{\mathrm{op}} % opposite category
\newcommand\id{\mathrm{id}} % identity morphism
\renewcommand\d{\mathrm{d}} % differential
\newcommand\comp{\fatsemi} % fat semicolon

\newcommand\Set{\mathbf{Set}} % category of sets
\newcommand\Euc{\mathbf{Euc}} % category of euclidean spaces
\newcommand\Mark{\mathbf{Mark}} % category of Markov kernels
\newcommand\Gauss{\mathbf{Gauss}} % category of Gaussian kernels
\newcommand\Conv{\mathbf{Conv}} % category of convex sets
\newcommand\Optic{\mathbf{Optic}} % category of optics

\title{Value Iteration is Optic Composition}
\def\titlerunning{Value Iteration is Optic Composition}
\author{Jules Hedges \and Riu Rodr\'iguez Sakamoto}
\def\authorrunning{Jules Hedges and Riu Rodr\'iguez Sakamoto}

\begin{document}

\maketitle

\begin{abstract}
	Dynamic programming is a class of algorithms used to compute optimal control policies for Markov decision processes. Dynamic programming is ubiquitous in control theory, and is also the foundation of reinforcement learning. In this paper, we show that value improvement, one of the main steps of dynamic programming, can be naturally seen as composition in a category of optics, and intuitively, the optimal value function is the limit of a chain of optic compositions. We illustrate this with three classic examples: the gridworld, the inverted pendulum and the savings problem. This is a first step towards a complete account of reinforcement learning in terms of parametrised optics.  \end{abstract}

\section{Introduction}

In this paper we describe basic concepts of dynamic programming in terms of categories of optics. The class of models we consider are discrete-time Markov decision processes, aka. discrete-time controlled Markov chains. There are classical methods of computing optimal control policies, underlying much of both classical control theory and modern reinforcement learning, known collectively as \emph{dynamic programming}. These are based on two operations that can be interleaved in many different ways: \emph{value improvement} and \emph{policy improvement}. The central idea of this paper is the slogan \emph{value improvement is optic precomposition}, or said differently, \emph{value improvement is a representable functor on optics}.

Given a control problem with state space $X$, a \emph{value function} $V : X \to \R$ represents an estimate of the long-run payoff of following a policy starting from any state, and can be equivalently represented as a costate $V : \binom{X}{\R} \to I$ in a category of optics. Every control policy $\pi$ also induces an optic $\lambda (\pi) : \binom{X}{\R} \to \binom{X}{\R}$. The general idea is that the forwards pass of the optic is a morphism $X \to X$ describing the dynamics of the Markov chain given the policy, and the backwards pass is a morphism $X \otimes \R \to \R$ which given the current state and the \emph{continuation payoff}, describing the total payoff from all future stages, returns the total payoff for the current stage given the policy, plus all future stages.

Given a policy $\pi$ and a value function $V : \binom{X}{\R} \to I$, the costate $\binom{X}{\R} \overset{\lambda (\pi)}\longrightarrow \binom{X}{\R} \overset{V}\longrightarrow I$ is a closer approximation of the value of $\pi$. This is called \emph{value improvement}. Iterating this operation 
\[ \ldots \binom{X}{\R} \overset{\lambda (\pi)}\longrightarrow \binom{X}{\R} \overset{\lambda (\pi)}\longrightarrow \binom{X}{\R} \overset{V}\longrightarrow I \]
converges efficiently to the true value function of the policy $\pi$.

Replacing $\pi$ with a new policy that is optimal for its value function is called \emph{policy improvement}. Repeating these steps is known as \emph{policy iteration}, and converges to the optimal policy and value function.

Alternatively, instead of repeating value improvement until convergence before each step of policy improvement, we can also alternate them, giving the composition of optics
\[ \ldots \binom{X}{\R} \overset{\lambda (\pi_2)}\longrightarrow \binom{X}{\R} \overset{\lambda (\pi_1)}\longrightarrow \binom{X}{\R} \overset{V}\longrightarrow I \]
where each policy $\pi_i$ is optimal for the value function to the right of it. This is known as \emph{value iteration}, and also converges to the optimal policy and value function.
For an account of convergence properties of these algorithms, classic textbooks are \cite[Sec.6]{Puterman}, \cite[Ch.1]{bertsekas_dpoc_vol2}.

In this paper we illustrate this idea, using mixed optics to account for the categorical structure of transitions in a Markov chain and the convex structure of expected payoffs, which typically form the kleisli and Eilenberg-Moore categories of a probability monad. This paper is partially intended as an introduction to dynamic programming for category theorists, focussing on illustrative examples rather than on heavy theory.

\subsection{Related work}

The precursor of this paper was early work on value iteration using open games \cite{hedges_etal_compositional_game_theory}. The idea originally arose around 2016 during discussions of the first author with Viktor Winschel and Philipp Zahn. An early version was planned as a section of \cite{hedges_morphisms_open_games} but cut partly for page limit reasons, and partly because the idea was quite uninteresting until it was understood how to model stochastic transitions in open games \cite{bolt_hedges_zahn_bayesian_open_games} via optics \cite{Riley}. In this paper we have chosen to present the idea without any explicit use of open games, both in order to clarify the essential idea and also to bring it closer to the more recent framework of categorical cybernetics \cite{towards_foundations_categorical_cybernetics}, which largely subsumes open games \cite{capucci_etal_translating_extensive_form}. (Although, actually using this framework properly is left for future work.)

A proof-of-concept implementation of value iteration with open games was done in 2019 by the first author and Wolfram Barfuss\footnote{Source currently available at \url{https://github.com/jules-hedges/open-games-hs/blob/og-v0.1/src/OpenGames/Examples/EcologicalPublicGood/EcologicalPublicGood.hs}}, implementing a model from \cite{barfuss_learning_dynamics} - a model of the social dilemma of emissions cuts and climate collapse as a stochastic game, or jointly controlled MDP - and verifying it against Barfuss' Matlab implementation. A far more advanced implementation of reinforcement learning using open games was developed recently by Philipp Zahn, currently closed-source, and was used for the paper \cite{eschenbaum_etal_robust_algorithmic_collusion}.

The most closely related work to ours is \cite{DJM}, which formulates MDPs in terms of F-lenses \cite{Spivak} of the functor $\operatorname{BiKl}(C\times -,\Delta(\RR\times -))^{\op}$, where $C\times-$ is the reader comonad and $\Delta(\RR\times-)$ is a probability monad over actions with their expected value.
A MDP there is a lens from states and potential state changes and rewards to the agents observation and input $\binom{X}{\Delta(X\times\R)}\to\binom{O}{I}$.
Our approach differs in two ways.
We firstly assume that the readout function is the identity, as we are not dealing with partial observability \cite{POMDP}.
Secondly, we specify a concrete structure of the backwards update map $f^*: X\times I\to \Delta(X\times R)$, which allows us to rearrange the interface of this lens from policies to value functions.
Doing so opens up the possibility of composing these lenses sequentially, which is the heart of the dynamic programming approach explored in this paper.

Bakirtzis et al. propose a category of MPDs as models of tasks \cite{bakirtzis_categorical_semantics_rl}. This emphasis on models allow them to compose different MPDs using fiber products and pushouts, and is agnostic to the control and RL algorithms that operate on them, which they take as given. Since our work focuses on a particular family of algorithms, we believe this approach is orthogonal to ours, and both could potentially be done simultaneously.

Another approach is to model MDPs as coalgebras from states to rewards and potential transitions, as done by Feys et al. \cite{LTVMDPs}.
They observe that the Bellman optimality condition for value iteration is a certain coalgebra-to-algebra morphism. We similarly believe this is orthogonal to our work and could potentially be unified.

A series of papers by Botta et al (for example \cite{botta13}) formulates dynamic programming in dependent type theory, accounting in a serious way for how different actions can be available in different states, a complication that we ignore in this paper. It may be possible to unify these approaches using dependent optics \cite{braithwaite_etal_fibre_optics,vertechi_dependent_optics}.

Finally, \cite{baez_erbele_categories_control} builds a category of signal flow diagrams, a widely used tool in control theory. Besides the common application to control theory there is little connection to this paper. In particular, time is implicit in their string diagrams, meaning their models have continuous time, whereas our approach is inherently discrete time. Said another way, composition in their category is `space-like' whereas ours is `time-like' - their morphisms are (open) systems whereas ours are processes.

\section{Dynamic programming}

\subsection{Markov Decision Processes}

A \emph{Markov decision process} (MDP) consists of a state space $X$, an action space $A$, a state transition function $f: X\times A\to X$, and a utility or reward function $U: X\times A\to \RR$. The state transition function is often taken to be stochastic, that is, to be given by probabilities $f (x' \mid x, a)$. In the stochastic case the utility function can be taken without loss of generality to be an expected utility function. We imagine actions to be chosen by an agent, who is trying to \emph{control} the Markov chain with the objective of optimising the long-run reward.

A \emph{policy} for an MDP is a function $\pi : X \to A$, which can also be taken to be either deterministic or stochastic. The type of policies encodes the Markov property: the choice of action depends only on the current state, and may not depend on any memory of past states.

Given an initial state $x_0 \in X$, a policy $\pi$ determines (possibly stochastically) a sequence of states
\[ x_0,\quad x_1=f(x_0,\pi(x_0)),\quad x_2=f(x_1,\pi(x_1)),\quad \dots \]
The total payoff is given by an infinite geometric sum of individual payoffs for each transition:
\begin{equation} \label{eq:policy_value}
    V_\pi(x_0) = \sum_{k=0}^\infty \beta^kU(x_k,\pi(x_k))
\end{equation}
where $0 < \beta < 1$ is a fixed \emph{discount factor} which balances the relevance of present and future payoffs. (There are other methods of obtaining a single objective from an infinite sequence of transitions, such as averaging, but we focus on discounting in this paper.)
A key idea behind dynamic programming is that this geometric sum can be equivalently written as a telescoping sum:
\[ V_\pi (x_0) = U (x_0, \pi (x_0)) + \beta (U (x_1, \pi (x_1)) + \beta (U (x_2, \pi (x_2) + \cdots))) \]

The \emph{control problem} is to choose a policy $\pi$ in order to maximise (the expected value of) $V_\pi (x_0)$.
In terms of decision theory, we assume that the agent choosing the policy operates under rational behaviour.
Continuous and independent preferences of outcome implies by the von Neumann-Morgenstern expected utility theorem that the utility function has as codomain the reals.

\subsection{Deterministic dynamic programming}

In dynamic programming, the agent's objective of maximizing the overall utility can be divided into two orthogonal goals:
to determine the value of a given policy $\pi$ (which we call the \emph{value improvement} step), and to determine the optimal policy $\pi^*$ (the \emph{policy improvement} step).
Bellman's equation is used as an update rule for both:
\begin{align}
    \text{\textbf{Value improvement:}} &&&V'(x) = U(x,\pi(x)) + \beta V(f(x,\pi(x))) \label{eq:value_improvement} \\
    \text{\textbf{Policy improvement:}} &&&\pi'(x) = \arg\max_{a \in A} U(x,a) + \beta V(f(x,a)) \label{eq:policy_improvement}
\end{align}
A Bellman optimality condition on the other hand determines the fixpoint of this update rule, and is met when $V'=V$ and $\pi'=\pi$ respectively.
% TODO: $\pi$ doesn't appear in Policy improvement. How can we describe an old value function of a new policy?

The update rule \eqref{eq:value_improvement} is the discounted sum \eqref{eq:policy_value} where the stream of states is co-recursively fixed by the policy $\pi$ and transition function $f$. % maybe op-recursive?
The co-recursive structure refers to the calculation of the utility of a state $x$, where one needs the utility of the \emph{next} state, while in a recursive structure, $x$ needs the \emph{previous} state, starting from an initial state as a base case.

Two classical algorithms use these two steps differently:
Policy iteration iterates value improvement until the current policy value is optimal before performing a policy improvement step, and value iteration interleaves both steps one after another.

% Finding an optimal policy $\pi^*$ such that $V_\pi(s_0)$ is maximized given a specific $s_0$ is always possible (local convergence). However, an optimal policy for any starting state (global convergence) might not exist \cite{Puterman}.

In policy iteration, a initial value function is chosen (usually $V(x)=0$), and a randomly chosen policy $\pi$ is evaluated by \eqref{eq:value_improvement} repeatedly until the value reaches a fixpoint, which is assured by the contraction mapping \cite{denardo_contraction_in_dp}.
% Check convergence formulation?
Once $V$ reaches (or in practice gets close to) a fixpoint $V'=V$ or another convergence condition, the policy improvement step \eqref{eq:policy_improvement} chooses a greedy policy as an improvement to $\pi$.
% The `policy iteration' refers to the repeated cycle of getting a new policy and evaluating it, until a new greedy improvement to the policy doesn't change it:
% (this is Bellman's optimality condition when the environment is probabilistic):
% \[ \pi_{j+1}(s) = \arg\max_a q_{\pi_j}(s,a) \]

A $q$-function or \emph{state-action value function} $q_\pi : X \times A \to \R$ describes the value of being in state $x$ and then taking action $a$, assuming that subsequent actions are taken by the policy $\pi$
\begin{equation}\label{eq:q}
    q_\pi(x,a) = U(x,a) + \beta V(f(x,a))
\end{equation}
The \emph{policy improvement theorem} \cite{bellman_1957} states that if a pair of deterministic policies $\pi,\pi': X\to A$ satisfies for all $x\in X$
\[ q_\pi(x,\pi'(x)) \geq V_\pi(x) \]
then $V_{\pi'}(x)\geq V_\pi(x)$ for all $x\in X$.

The optimal policy $\pi^*$, if it exists, is the policy which if followed from any state, generates the maximum value.
This is a Bellman optimality condition which fuses the two steps \eqref{eq:value_improvement}, \eqref{eq:policy_improvement}:
\begin{equation} \label{eq:Bellman}
    V_{\pi^*}(x) = \max_{a \in A} U(x,a) + \beta V_{\pi^*}(f(x,a))
\end{equation}

\emph{Value iteration} is a special policy iteration algorithm insofar it stops the update rule for value improvement to one step, by truncating the sum \eqref{eq:policy_value} to the first summand.
Moreover, it introduces the value improvement step implicitly in the policy improvement, which assigns a value to states
\[ V'(x) = \max_{a \in A} U(x,a) + \beta V(f(x,a)) \]
while the policy in each iteration is still recoverable as
\[ \pi'(x) = \arg\max_{a \in A} U(x,a) + \beta V(f(x,a)) \]

\subsection{Stochastic dynamic programming}

Stochasticity can be introduced in different places in a MDP:
\begin{enumerate}%[label=(\alph*)]
    \item in the policy $\pi: X\to \Delta A$, where the probability of the policy $\pi$ taking action $a$ in a state $x$ is now notated $\pi(a \mid x)$.
    \item in the transition function $f: X\times A\to \Delta X$ and potentially the reward function $U: X\times A\to \Delta\R$ independently.
    \item usually the reward is included inside the transition function $f: X\times A\to \Delta(X\times\R)$, allowing correlated next states and rewards.
    This is relevant when the reward is morally from the next state, rather than the current state and action.
    If the reward were truly from the current state and action, the transition function can be decomposed into a function $f: X\times A\to \Delta X \times \Delta \R$.
\end{enumerate}
In this section we assume for simplicity that $\Delta$ is the finite support distribution monad, although the equations in the following can be formulated for arbitrary distributions by replacing the sum with an appropriate integral.
% Alternative explanation: the MDP sits inside a Markov category, in which the transition and utility functions are not deterministic morphisms: they don't commute with the copy map (comultiplication) \cite[Def.10.1]{Synthetic_approach}

The policy value update rule \eqref{eq:value_improvement} becomes stochastic, and adopts a slightly different form depending on which part of the MDP is stochastic. For the cases 1. and 2.:
\begin{align}
    V'(x) &= \sum_a\pi(a \mid x)(U(x,a) + \beta V(f(x,a))) \\
    V'(x) &= \sum_r U(r \mid x,a)r + \sum_{x'}f(x' \mid x,a)\beta V(x') \\
    \intertext{In the most general case, that is 1. together with 3.:}
    V'(x) &= \sum_{a\in A}\pi(a \mid x)\sum_{x',r}f(x',r \mid x,a)(r + \beta V(x')) \label{eq:prediction}
\end{align}
(Note that the sum over $r$ is over the support of $f (- \mid x, a)$, which we assume here to be finite, although in general it can be replaced with an integral.)

%Note that even though that $\Delta$ is a finite-support probability monad, it can still be applied to $\RR$. $\Delta\RR$ restricts the subset of all possible distributions over $\RR$ to just the set of finite-support ones, and therefore we can sum over values of $r'$ as there are only finitely many.

% Expectation values appear in Bellman's equation when the transition function and/or the utility function are probabilistic and we need to take the expected transition: $\mathbb{E}_\pi$ means, given a fixed policy $\pi$, take expectations over everything else.
% \[ v_\pi^{k+1}(s)=\sum_{s',r}p(s',r|s,\pi(s))[r+\beta v_\pi^k(s')] \]

% \subsubsection{Policy improvement}
The policy improvement theorem holds in the stochastic setting \cite[Sec.4.2]{RL_Intro} by defining
\[ q_\pi(s,\pi'(s))=\sum_a \pi'(a \mid s)q_\pi(s,a) \]

\subsection{Gridworld example}

A classic example in reinforcement learning is the Gridworld environment, where an agent moves in the four cardinal directions in a rectangular grid.
States of this finite MDP correspond to the positions that the agent can be in.

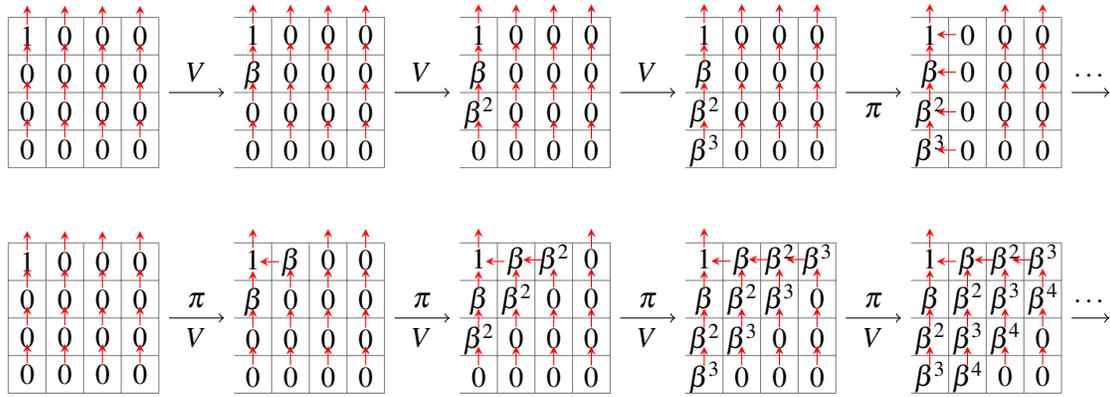
\begin{figure}[ht!]
    \centering
    \begin{tikzpicture}
        \tikzstyle{arrow} = [red,-stealth,shorten <=-3pt,shorten >=-3pt]
        \tikzstyle{matrixsty} = [matrix of nodes,nodes={inner sep=0pt,text width=.5cm,align=center,minimum height=.5cm}]
        \matrix(m0)[matrixsty] at (1,1){
        1 & 0 & 0 & 0 \\
        0 & 0 & 0 & 0 \\
        0 & 0 & 0 & 0 \\
        0 & 0 & 0 & 0 \\};
            
        \matrix(m1)[matrixsty] at (4,1){
        1 & 0 & 0 & 0 \\
        $\beta$ & 0 & 0 & 0 \\
        0 & 0 & 0 & 0 \\
        0 & 0 & 0 & 0 \\};

        \matrix(m2)[matrixsty] at (7,1){
        1 & 0 & 0 & 0 \\
        $\beta$ & 0 & 0 & 0 \\
        $\beta^2$ & 0 & 0 & 0 \\
        0 & 0 & 0 & 0 \\};

        \matrix(m3)[matrixsty] at (10,1){
        1 & 0 & 0 & 0 \\
        $\beta$ & 0 & 0 & 0 \\
        $\beta^2$ & 0 & 0 & 0 \\
        $\beta^3$ & 0 & 0 & 0 \\};

        \matrix(m4)[matrixsty] at (13,1){
        1 & 0 & 0 & 0 \\
        $\beta$ & 0 & 0 & 0 \\
        $\beta^2$ & 0 & 0 & 0 \\
        $\beta^3$ & 0 & 0 & 0 \\};

        % Second row: Value Iteration
        \matrix(m5)[matrixsty] at (1,-2){
        1 & 0 & 0 & 0 \\
        0 & 0 & 0 & 0 \\
        0 & 0 & 0 & 0 \\
        0 & 0 & 0 & 0 \\};
            
        \matrix(m6)[matrixsty] at (4,-2){
        1 & $\beta$ & 0 & 0 \\
        $\beta$ & 0 & 0 & 0 \\
        0 & 0 & 0 & 0 \\
        0 & 0 & 0 & 0 \\};

        \matrix(m7)[matrixsty] at (7,-2){
        1 & $\beta$ & $\beta^2$ & 0 \\
        $\beta$ & $\beta^2$ & 0 & 0 \\
        $\beta^2$ & 0 & 0 & 0 \\
        0 & 0 & 0 & 0 \\};

        \matrix(m8)[matrixsty] at (10,-2){
        1 & $\beta$ & $\beta^2$ & $\beta^3$ \\
        $\beta$ & $\beta^2$ & $\beta^3$ & 0 \\
        $\beta^2$ & $\beta^3$ & 0 & 0 \\
        $\beta^3$ & 0 & 0 & 0 \\};

        \matrix(m9)[matrixsty] at (13,-2){
        1 & $\beta$ & $\beta^2$ & $\beta^3$ \\
        $\beta$ & $\beta^2$ & $\beta^3$ & $\beta^4$ \\
        $\beta^2$ & $\beta^3$ & $\beta^4$ & 0 \\
        $\beta^3$ & $\beta^4$ & 0 & 0 \\};

        % grids
        \foreach \gridx in {0,3,...,12} {
            \pgfmathtruncatemacro{\gridxx}{\gridx+2}
            \draw[step=.5cm,color=gray] (\gridx,0) grid (\gridxx,2);
            \draw[step=.5cm,color=gray] (\gridx,-3) grid (\gridxx,-1);
        }
        
        % policy arrows for first four grids of Policy Iteration
        \foreach \gridid in {m0,m1,m2,m3} {
            \foreach \row in {1,...,4} {
                \foreach \col in {1,...,4} {
                    \draw[arrow] (\gridid-\row-\col) -- ++ (0,.3cm);
                }
            }
        }
        % policy arrows for last grid of Policy Iteration (m4)
        \foreach \col in {1,3,4} {
            \foreach \row in {1,...,4} {
                \draw[arrow] (m4-\row-\col) -- ++ (0,.3cm);
            }
        }
        \foreach \row in {1,...,4} {
            \draw[arrow] (m4-\row-2) -- ++ (-.3cm,0);
        }

        % policy arrows for all grids of Value Iteration, all rows except first
        \foreach \gridid in {m5,m6,m7,m8,m9} {
            \foreach \row in {2,3,4} {
                \foreach \col in {1,...,4} {
                    \draw[arrow] (\gridid-\row-\col) -- ++ (0,.3cm);
                }
            }
        }
        % upwards policy arrows for first row in all grids of VI
        \foreach \gridid/\col in {m5/1,m5/2,m5/3,m5/4,m6/1,m6/3,m6/4,m7/1,m7/4,m8/1,m9/1} {
            \draw[arrow] (\gridid-1-\col) -- ++ (0,.3cm);
        }
        % left facing policy arrow for first row in all grids of VI
        \foreach \gridid/\col in {m6/2,m7/2,m7/3,m8/2,m8/3,m8/4,m9/2,m9/3,m9/4} {
            \draw[arrow] (\gridid-1-\col) -- ++ (-.3cm,0);
        }

        % arrows between grids of Policy Iteration
        \draw[->] (m0.east) -- (m1.west) node[above,midway]{$V$};
        \draw[->] (m1.east) -- (m2.west) node[above,midway]{$V$};
        \draw[->] (m2.east) -- (m3.west) node[above,midway]{$V$};
        \draw[->] (m3.east) -- (m4.west) node[below,midway]{$\pi$};
        \draw[->] (m4.east) -- ++ (.5cm,0) node[above,midway]{$\cdots$};
        % arrows between grids of Value Iteration
        \draw[->] (m5.east) -- (m6.west) node[above,midway]{$\pi$} node[below,midway]{$V$};
        \draw[->] (m6.east) -- (m7.west) node[above,midway]{$\pi$} node[below,midway]{$V$};
        \draw[->] (m7.east) -- (m8.west) node[above,midway]{$\pi$} node[below,midway]{$V$};
        \draw[->] (m8.east) -- (m9.west) node[above,midway]{$\pi$} node[below,midway]{$V$};
        \draw[->] (m9.east) -- ++ (.5cm,0) node[above,midway]{$\cdots$};

    \end{tikzpicture}
    \caption{Difference between policy iteration (above) and value iteration (below). The numbers in the cells are state values and the red arrows are the directions dictated by the policy at each stage. The arrows between grids indicate what kind of update the algorithm does, either value improvement ($V$) or policy improvement ($\pi$). Notice how policy iteration performs value improvement three times before updating the policy, whereas value iteration improves the value and the policy at each stage.}
    \label{fig:Gridworld}
\end{figure}

Assume that all transitions and policies are deterministic, and that the transition function prevents the agent from moving outside the boundary.
Suppose that the environment rewards 0 value for all states except the top left corner, where the reward is 1 (see figure \ref{fig:Gridworld}).

Starting with a policy which moves upwards in all states and a value function which rewards 1 only in the top left corner, a policy iteration algorithm would improve the value of the current policy until converging to the optimal values in the leftmost column, before updating the policy, while a value iteration algorithm would update the value function and also update the policy.

Take the finite set of positions as the state space $X$, and $A=\{\leftarrow,\rightarrow,\uparrow,\downarrow\}$ as the action space.

This example can be made stochastic if we add stochastic policies like $\epsilon$-greedy, where the action that the agent takes is the one with maximum value with probability $1 - \epsilon$ and a random one with probability $\epsilon$.
Another way is for the transition function to be stochastic, for example with a wind current that shifts the next state to the right with some probability $\epsilon$.

\subsection{Inverted pendulum example}

A task that illustrates a continuous state space MDP is the control of a pendulum balanced over a cart, which can be described in continuous-time exactly by two non-linear differential equations \cite[Example 2E]{control_systems_design}:
\begin{align*}
    (M+m)\ddot{y}+mL\ddot{\theta}\cos\theta-mL\dot{\theta}^2\sin\theta &= a \\
    mL\ddot{y}\cos\theta+mL^2\ddot{\theta}-mLg\sin\theta &= 0
\end{align*}
where $M$ is the mass of the cart, $m$ the mass of the pendulum, $L$ the length of the pendulum, $\theta$ the angle of the pendulum with respect to the upwards position, $y$ the carts horizontal position, $g$ the gravitational constant and $a$ is our control function (usually denoted $u$).
We rewrite the state variables as $x=[y,\dot{y},\theta,\dot{\theta}]^\top$.

Sampling the trajectory of continuous-time dynamics $\frac{\d}{\d t}x(t)=f(x(t))$ by $x_k=x(k\Delta t)$, one can define the discrete-time propagator $F_{\Delta t}$ by
\[ F_{\Delta t}(x(t))= x(t) + \int_t^{t+\Delta t} f(x(\tau)) \d\tau \]
which allows to model the system with $x_{k+1}=F_{\Delta t} (x_k)$.

A more common approach is to observe that the system of equations $\dot{x}=A(x)+B(x)a$ with $A$ and $B$ being non-linear functions of the state space, can be \emph{linearized} near a (not necessarily stable) equilibrium state, like the pendulum being in the upwards position.
There we can assume certain approximations like $\cos \theta\approx 1$ and $\sin\theta\approx \theta$, as well as small velocities leading to negligible quadratic terms $\dot{\theta}^2\approx 0$ and $\dot{y}^2\approx 0$.
This linearization around a fixpoint allows for the expression $\dot{x}=A x(t)+B a(t)$, where the matrix $A$ and vector $B$ are constants given by
\begin{align*}
    A &= \begin{pmatrix}
        0 & 1 & 0 & 0 \\
        0 & 0 & -\frac{mg}{M} & 0 \\
        0 & 0 & 0 & 1 \\
        0 & 0 & \frac{(M+m)g}{ML} & 0
    \end{pmatrix} & B &= \begin{pmatrix}
        0 \\
        \frac{1}{M} \\
        0 \\
        -\frac{1}{ML}
    \end{pmatrix}
\end{align*}

If we assume that the observation of the pendulum angle and cart position is discretized in time, an a priori time-discretization of this model using Euler approximation follows $x_{k+1}=x_k+\Delta t (A x_k+ B a_k)$, with the same constants, where $k$ indexes time steps.
Therefore we can say that the time-discretized, linearized model of the inverted pendulum over a cart follows a deterministic MDP for which a controller $u$ can be learned.
We take the state space as $\R^4$ and the action space of the force exerted to the cart as $\R$.

The time-discretised formulation of this problem is more common in reinforcement learning settings than in `classical' control theory. In that case, a a common payoff function is to obtain one unit of reward for each time step that the pendulum is maintained within a threshold of angles. The not linearized, not time-discretized setting, which is more common in optimal control theory, allows the reward, which is usually termed negatively as a cost function $J$, to have a much more flexible expression, in terms of time spent towards the equilibrium, energy spent to control the device, etc.
\[ J(x,a) = \int_0^\infty C(x(t),a(t)) \d t \]

\subsection{Savings problem example}

The \emph{savings problem} is one of the most important models in economics, modelling the dilemma between saving and consumption of resources \cite[part IV]{ljungqvist_sargent_recursive_macroeconomic_theory} (see also \cite{stokey_lucas_recursive_methods_economic_dynamics}). It is also mathematically closely related to the problem of charging a battery, for example choosing when to draw electricity from a power grid to raise the water level in a reservoir \cite{mody_steffen_optimal_charging_electric_vehicles}. 

At each discrete time step $k$, an agent receives an income $i_k$. They also have a bank balance $x_k$, which accumulates interest over time (this could also be, for example, an investment portfolio yielding returns). At each time step the agent makes a choice of \emph{consumption}, which means converting their income into utility (or, more literally, things from which they derive utility). If the consumption in some stage is less than their income then the difference is added to the bank balance, and if it is more than the difference is taken from the bank balance. The dilemma is that the agent receives utility only from consumption, but saving gives the possibility of higher consumption later due to interest. The optimal balance between consumption and saving depends on the discount factor, which models the agent's preference between consumption now and consumption in the future.

In the most basic version of the model, all values can be taken as deterministic, and the income $i_k$ can also be taken as constant. This basic model can be expanded in many ways, for example with forecasts and uncertainty about income and interest rates. A straightforward extension, which we will consider in this paper, is that income is normally distributed $i \sim \mathcal N (\mu, \sigma)$, independently in each time step.

We take the state space and action space both as $X = A = [0, \infty)$. Given the current bank balance $x$ and consumption decision $a$, the utility in the current stage is $U (x, a) = \min \{ a, x + i \}$. (That is, the agent's consumption is capped by their current bank balance.) The state transition is given by $f (x, a) = \max \{ (1 + \gamma) x - a + i, 0 \}$, where $\gamma$ is the interest rate.

\section{Optics}

In this section we recall material on categories of mixed optics, mostly taken from \cite{roman_etal_profunctor_optics}.

\subsection{Categories of optics}

Given a monoidal category $\M$ and a category $\C$, an action of $\M$ on $\C$ is given by a functor $\bullet : \M \times \C \to \C$ with coherence isomorphisms $I \bullet X \cong X$ and $(M \otimes N) \bullet X \cong M \bullet (N \bullet X)$. $\C$ is called an $\M$-actegory.

Given a pair of $\M$-actegories $\C,\D$, we can form the category of optics $\Optic_{\C,\D}$. Its objects are pairs $\binom{X}{X'}$ where $X$ is an object of $\C$ and $X'$ is an object of $\D$. Hom-sets are defined by the coend
\[ \Optic_{\C,\D} \left( \binom{X}{X'}, \binom{Y}{Y'} \right) = \int^{M : \M} \C (X, M \bullet Y) \times \D (M \bullet Y', X') \]
in the category $\Set$. Such a morphism is called an optic, and consists of an equivalence class of triples $(M, f, f')$ where $M$ is an object of $\M$, $f : X \to M \bullet Y$ in $\C$ and $g : M \bullet Y' \to X'$ in $\D$. We call $M$ the residual, $f$ the forwards pass and $f'$ the backwards pass, so we think of the residual as mediating communication from the forward pass to the backward pass. Composition of optics works by taking the monoidal product in $\M$ of the residuals.

A common example is a monoidal category $\M = \C = \D$ acting on itself by the monoidal product, so
\[ \Optic_\C \left( \binom{X}{X'}, \binom{Y}{Y} \right) = \int^{M : \C} \C (X, M \otimes Y) \times \C (M \otimes Y', X') \]
This is the original definition of optics from \cite{Riley}. If $\C$ is additionally cartesian monoidal then we can eliminate the coend to produce \emph{concrete lenses}:
\begin{align*}
	\int^{M : \C} \C (X, M \times Y) \times \C (M \times Y', X') &\cong \int^{M : \C} \C (X, M) \times \C (X, Y) \times \C (M \times Y', X') \\
	&\cong \C (X, Y) \times \C (X \times Y', X')
\end{align*}
On the other hand, if $\C$ is monoidal closed then we can eliminate the coend in a different way to produce \emph{linear lenses}:
\begin{align*}
	\int^{M : \C} \C (X, M \otimes Y) \times \C (M \otimes Y', X') &\cong \int^{M : \C} \C (X, M \otimes Y) \times \C (M, [Y', X']) \\
	&\cong \C (X, [Y', X'] \otimes Y)
\end{align*}
Both of these proofs use the \emph{ninja Yoneda lemma} for coends \cite{loregian-coend-cofriend}.

\begin{example}
	Let $\Set$ act on itself by cartesian product. Optics $\binom{X}{X'} \to \binom{Y}{Y'}$ in $\Optic_\Set$ can be written equivalently as pairs of functions $X \to Y$ and $X \times Y' \to X'$, or as a single function $X \to Y \times (Y' \to X')$.
\end{example}

\begin{example}
	Let $\Euc$ be the category of Euclidean spaces and smooth functions, which is cartesian but not cartesian closed. Optics $\binom{X}{X'} \to \binom{Y}{Y'}$ in $\Optic_\Euc$ can be written as pairs of smooth functions $X \to Y$ and $X \times Y' \to X'$.
\end{example}

\begin{example}
	Let $\Mark$ be the category of sets and finite support Markov kernels, which is the kleisli category of the finite support probability monad $\Delta : \Set \to \Set$. It is a prototypical example of a Markov category \cite{Synthetic_approach}, and it is neither cartesian monoidal nor monoidal closed. Optics $\binom{X}{X'} \to \binom{Y}{Y'}$ in $\Optic_\Mark$ can only be written as optics, it is not possible to eliminate the coend. This is the setting used for Bayesian open games \cite{bolt_hedges_zahn_bayesian_open_games}.
\end{example}

\begin{example}
	Let $\Conv$ be the category of convex sets, which is the Eilenberg-Moore category of the finite support probability monad \cite{fritz09}. A convex set can be thought of a set with an abstract expectation operator $\mathbb E : \Delta X \to X$. Thus the functor $\Delta : \Mark \to \Conv$ given by $X \mapsto \Delta (X)$ on objects is fully faithful. $\Conv$ has finite products which are given by tupling in the usual way. $\Conv$ also has a closed structure: the set of convex functions $X \to Y$ themselves form a convex set $[X, Y]$ pointwise. However $\Conv$ is not cartesian closed: instead there is a different monoidal product making it monoidal closed \cite[section 2.2]{sturtz_categorical_probability} (see also \cite{kock_closed_categories_commutative_monads}). This monoidal product ``classifies biconvex maps'' in the same sense that the tensor product of vector spaces classifies bilinear maps. The embedding $\Delta : \Mark \to \Conv$ is strong monoidal for this monoidal product, not for the cartesian product of convex sets.
	
	We can define an action of $\Mark$ on $\Conv$, by $M \bullet X = \Delta (M) \otimes X$ \cite[section 5.5]{capucci_gavranovic_actegories}. Together with the self-action of $\Mark$, we get a category $\Optic_{\Mark,\Conv}$ given by
	\begin{align*}
		\Optic_{\Mark, \Conv} \left( \binom{X}{X'}, \binom{Y}{Y'} \right) &= \int^{M : \Mark} \Mark (X, M \otimes Y) \times \Conv (\Delta (M) \otimes Y', X') \\
		&\cong \int^{M : \Mark} \Mark (X, M \otimes Y) \times \Conv (\Delta (M), [Y', X'])
	\end{align*}
	(This coend cannot be eliminated because the embedding $\Delta : \Mark \to \Conv$ does not have a right adjoint.)
	
	This category of optics will be very useful for Markov decision processes, where the forwards direction is a Markov kernel and the backwards direction is a function involving expectations.
\end{example}

\subsection{Monoidal structure of optics}

A category of optics $\Optic_{\C,\D}$ is itself (symmetric) monoidal, when $\C$ and $\D$ are (symmetric) monoidal in a way that is compatible with the actions of $\M$. The details of this have been recently worked out in \cite{capucci_gavranovic_actegories}. The monoidal product on objects of $\Optic_{\C,\D}$ is given by pairwise monoidal product. All of the above examples are symmetric monoidal.

A monoidal category of optics comes equipped with a string diagram syntax \cite{hedges_coherence_lenses_open_games}. This has directed arrows representing the forwards and backwards passes, and right-to-left bending wires but not left-to-right bending wires. The residual of the denoted optic can be read off from a diagram, as the monoidal product of the wire labels of all right-to-left bending wires. 
For example, a typical optic $(M, f, f') \in \Optic_\C \left( \binom{X}{X'}, \binom{Y}{Y'} \right)$ is denoted by the diagram
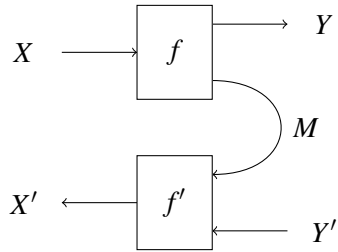
\begin{center}
	\begin{tikzpicture}
	\begin{pgfonlayer}{nodelayer}
		\node [style=2morph] (0) at (0, 2) {$f$};
		\node [style=2morph] (1) at (0, -2) {$f'$};
		\node [style=none] (2) at (-1, 2) {};
		\node [style=none] (3) at (-1, -2) {};
		\node [style=none] (4) at (1, 2.75) {};
		\node [style=none] (5) at (1, 1.25) {};
		\node [style=none] (6) at (1, -1.25) {};
		\node [style=none] (7) at (1, -2.75) {};
		\node [style=none] (8) at (-3, 2) {};
		\node [style=none] (9) at (-3, -2) {};
		\node [style=none] (10) at (3, 2.75) {};
		\node [style=none] (11) at (3, 1.25) {};
		\node [style=none] (12) at (3, -1.25) {};
		\node [style=none] (13) at (3, -2.75) {};
		\node [style=none] (14) at (-4, 2) {$X$};
		\node [style=none] (15) at (-4, -2) {$X'$};
		\node [style=none] (16) at (4, 2.75) {$Y$};
		\node [style=none] (17) at (4, -2.75) {$Y'$};
		\node [style=none] (18) at (3.5, 0) {$M$};
	\end{pgfonlayer}
	\begin{pgfonlayer}{edgelayer}
		\draw [style=arrow] (8.center) to (2.center);
		\draw [style=arrow] (3.center) to (9.center);
		\draw [style=arrow] (4.center) to (10.center);
		\draw [style=arrow] (13.center) to (7.center);
		\draw [style=arrow, in=0, out=0, looseness=2.50] (5.center) to (6.center);
	\end{pgfonlayer}
\end{tikzpicture}

\end{center}

These diagrams have only been properly formalised for a monoidal category acting on itself, so for mixed optics we need to be very careful and are technically being informal.

Costates in monoidal categories of optics, that is optics $\binom{X}{X'} \to I$ (where $I = \binom{I}{I}$ is the monoidal unit of $\Optic_{\C,\D}$), are a central theme of this paper. When we have a monoidal category acting on itself, costates in $\Optic_\C$ are given by
\[ \Optic_\C \left( \binom{X}{X'}, I \right) = \int^{M : \C} \C (X, M \otimes I) \times \C (M \otimes I, X') \cong \C (X, X') \]
Thus \emph{costates in optics are functions}. A different way of phrasing this is by defining a functor $K : \Optic_\C^\op \to \Set$ given on objects by $K \binom{X}{X'} = \C (X, X')$, and then showing that $K$ is representable \cite{hedges_morphisms_open_games}. We will generally treat this isomorphism as implicit, sometimes referring to costates as though they are functions.

In the case of a cartesian monoidal category $\C$, given a concrete lens $f : X \to Y$, $f' : X \times Y' \to X'$ and a function $k : Y \to Y'$, the action of $K$ gives us the function $X \to X'$ given by $x \mapsto f' (x, k (f (x)))$.

When we have $\M = \C$ acting on both itself and on $\D$ (which includes all of the examples above) then similarly
\[ \Optic_{\C,D} \left( \binom{X}{X'}, I \right) = \int^{M : \C} \C (X, M \otimes I) \times \D (M \bullet I, X') \cong \D (X \bullet I, X') \]

\section{Dynamic programming with optics}

Given an MDP with state space $X$ and action space $A$, we can convert it to an optic $\binom{X \otimes A}{\R} \to \binom{X}{\R}$. The category of optics in which this lives can be `customised' to some extent, and depends on the class of MDPs that we are considering and how much typing information we choose to include. The definition of this optic is given by the following string diagram:
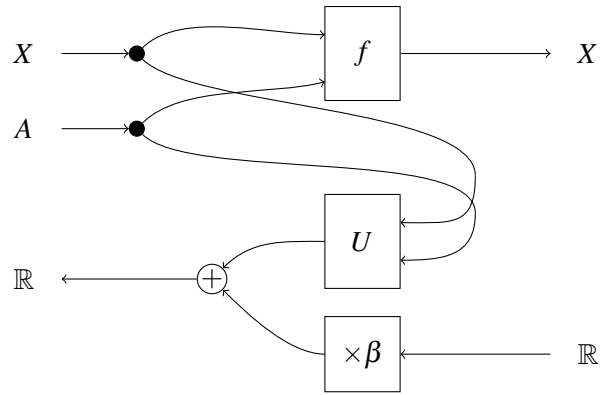
\begin{center}
	\begin{tikzpicture}
	\begin{pgfonlayer}{nodelayer}
		\node [style=copy2] (0) at (-2, 2) {};
		\node [style=copy2] (1) at (-2, 0) {};
		\node [style=2morph] (2) at (4, -3) {$U$};
		\node [style=copy] (3) at (0, -4) {$+$};
		\node [style=1morph1] (4) at (4, -6) {$\times \beta$};
		\node [style=2morph] (5) at (4, 2) {$f$};
		\node [style=none] (6) at (3, 2.5) {};
		\node [style=none] (7) at (3, 1.25) {};
		\node [style=none] (8) at (5, 2) {};
		\node [style=none] (9) at (9, 2) {};
		\node [style=none] (10) at (-4, 2) {};
		\node [style=none] (11) at (-4, 0) {};
		\node [style=none] (12) at (7, -1.25) {};
		\node [style=none] (13) at (7, -2.25) {};
		\node [style=none] (14) at (5, -2.5) {};
		\node [style=none] (15) at (5, -3.5) {};
		\node [style=none] (16) at (3, -3) {};
		\node [style=none] (17) at (3, -6) {};
		\node [style=none] (18) at (5, -6) {};
		\node [style=none] (19) at (9, -6) {};
		\node [style=none] (20) at (-4, -4) {};
		\node [style=none] (21) at (-5, 2) {$X$};
		\node [style=none] (22) at (-5, 0) {$A$};
		\node [style=none] (23) at (-5, -4) {$\mathbb R$};
		\node [style=none] (24) at (10, 2) {$X$};
		\node [style=none] (25) at (10, -6) {$\mathbb R$};
	\end{pgfonlayer}
	\begin{pgfonlayer}{edgelayer}
		\draw [style=arrow, in=180, out=45, looseness=0.75] (0) to (6.center);
		\draw [style=arrow, in=195, out=45, looseness=0.75] (1) to (7.center);
		\draw [style=arrow] (8.center) to (9.center);
		\draw [style=arrow] (10.center) to (0);
		\draw [style=arrow] (11.center) to (1);
		\draw [in=90, out=-45, looseness=0.50] (0) to (12.center);
		\draw [style=arrow, in=0, out=-90, looseness=1.50] (12.center) to (14.center);
		\draw [in=90, out=-45, looseness=0.50] (1) to (13.center);
		\draw [style=arrow, in=0, out=-90, looseness=1.25] (13.center) to (15.center);
		\draw [style=arrow, in=45, out=180] (16.center) to (3);
		\draw [style=arrow] (19.center) to (18.center);
		\draw [style=arrow, in=-45, out=180, looseness=0.75] (17.center) to (3);
		\draw [style=arrow] (3) to (20.center);
	\end{pgfonlayer}
\end{tikzpicture}

\end{center}
To be clear, this diagram is not completely formal because we are making some assumptions about the category of optics we work in. In general, we require the forwards category $\C$ to be a Markov category (giving us copy morphisms $\Delta_X$ and $\Delta_A$), and the backwards category $\D$ must have a suitable object $\R$ together with morphisms $\times \beta : \R \to \R$ and $+ : \R \otimes \R \to \R$. Specific examples of interpretations of this diagram will be explored below. When the forwards category acts on the backwards category, then the forwards pass is a morphism $g : X \otimes A \to X \otimes X \otimes A$ in $\C$ where
\[ g = \Delta_{X \otimes A}\comp (f \otimes \id_{X \otimes A}) \]
and the backwards pass is a morphism $g' : X \bullet A \bullet \R \to \R$ in $\D$ encoding the function $g' (x, a, r) = \mathbb E U (x, a) + \beta r$.
The resulting optic is given by $\lambda = (X \otimes A, g, g') : \binom{X \otimes A}{\R} \to \binom{X}{\R}$ in $\Optic_{\C, \D}$.

Given a policy $\pi : X \to A$, we lift it to an optic $\overline\pi : \binom{X}{\R} \to \binom{X \otimes A}{\R}$, by
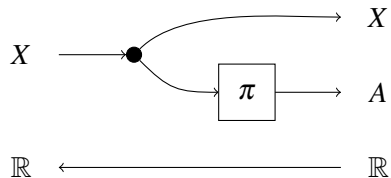
\begin{center}
	\begin{tikzpicture}
	\begin{pgfonlayer}{nodelayer}
		\node [style=none] (8) at (-3, 2) {};
		\node [style=none] (9) at (-3, -1) {};
		\node [style=none] (13) at (4.5, -1) {};
		\node [style=none] (14) at (-4, 2) {$X$};
		\node [style=none] (15) at (-4, -1) {$\R$};
		\node [style=none] (17) at (5.5, -1) {$\R$};
		\node [style=copy2] (18) at (-1, 2) {};
		\node [style=none] (21) at (2.75, 1) {};
		\node [style=none] (22) at (4.5, 3) {};
		\node [style=none] (23) at (4.5, 1) {};
		\node [style=none] (24) at (1.25, 1) {};
		\node [style=1morph] (25) at (2, 1) {$\pi$};
		\node [style=none] (26) at (5.5, 1) {$A$};
		\node [style=none] (27) at (5.5, 3) {$X$};
	\end{pgfonlayer}
	\begin{pgfonlayer}{edgelayer}
		\draw [style=arrow] (13.center) to (9.center);
		\draw [style=arrow, in=-180, out=45, looseness=0.75] (18) to (22.center);
		\draw [style=arrow] (21.center) to (23.center);
		\draw [style=arrow, in=-180, out=-45, looseness=1.25] (18) to (24.center);
		\draw [style=arrow] (8.center) to (18);
	\end{pgfonlayer}
\end{tikzpicture}

\end{center}
Here we are also assuming that the forwards category has a copy morphisms $\Delta_X$ (for example, because it is a Markov category), and the backwards category has a suitable object of real numbers. The interpretation of this diagram is the optic $(I, \Delta_X\comp (\id_X \otimes \pi), \id_\R)$.

%We define dynamic programming as a dynamical system on pairs consisting of a policy and a value function. Value improvement and policy improvement update each of these using the other. Fixpoints of this dynamical system characterise solutions of a Bellman equation. The fixpoint corresponds to an infinite composition of the backwards pass, corresponding to the telescoping sum
%\[ U (x_0, a_0) + \beta (U (x_1, a_1) + \beta (U (x_2, a_2) + \cdots)) \]
%which converges to the geometric sum $\sum_{i = 0}^\infty \beta^i U (x_i, a_i)$, where $x_i$ is the sequence of states and $a_i = \pi_i (x_i)$ is the sequence of actions.

\subsection{Discrete-space deterministic decision processes}\label{sec:4.1}

Consider a deterministic decision process with a discrete set of states $X$, discrete and finite set of actions $A$, transition function $f : X \times A \to X$, payoff function $U : X \times A \to \R$ and discount factor $\beta \in (0, 1)$. We convert this into an optic $\lambda = (X \times A, g, g') : \binom{X \times A}{\R} \to \binom{X}{\R}$ in $\Optic_\Set$, whose forwards pass is $g (x, a) = (x, a, f (x, a))$ and whose backwards pass is $g' (x, a, r) = U (x, a) + \beta r$.

Consider a dynamical system with the state space $A^X \times \Optic_\Set \left( \binom{X}{\R}, I \right)$. Elements of this are pairs $(\pi, V)$ of a policy $\pi : X \to A$ and a value function $V : X \to \R$. We can define two update steps:
\begin{align*}
	\text{\textbf{Value improvement:}}& &&(\pi, V) \mapsto (\pi, \overline\pi\comp \lambda\comp V) \\
	\text{\textbf{Policy improvement:}}& &&(\pi, V) \mapsto (x \mapsto \arg\max_{a \in A} (\lambda\comp V) (x, a), V)
\end{align*}
(We assume that $\arg\max$ is canonically defined, for example because $A$ is equipped with an enumeration so that we can always choose the first maximiser.)
	
Unpacking and applying the isomorphism between costates in lenses and functions, a step of value improvement replaces $V$ with
\[ V' (x) = U (x, \pi (x)) + \beta V (f (x, \pi (x))) \]
and a step of policy improvement replaces $\pi$ with
\[ \pi' (x) = \arg\max_{a \in A} U (x, a) + \beta V (f (x, a)) \]
Iterating the value improvement step converges to a value function which is the optimal value function for the current (not necessarily optimal) policy $\pi$.
A fixpoint of alternating steps of value improvement and policy improvement is a pair $(\pi^*, V^*)$ satisfying
\begin{align*}
	V^* (x) &= \max_{a \in A} (\lambda\comp V^*) (x, a) = \max_{a \in A} U (x, a) + \beta V^* (f (x, a)) \\
	\pi^* (x) &= \arg\max_{a \in A} (\lambda\comp V^*) (x, a) = \arg\max_{a \in A} U (x, a) + \beta V^* (f (x, a))
\end{align*}
%The first of these two equations is a \emph{Bellman equation}.

\begin{example}
	[Gridworld example]
    A policy of an agent in our version of Gridworld (Figure \ref{fig:Gridworld}) is a function from the $4\times 4$ set of states $X = \{ 1, 2, 3, 4 \}^2$ that we index by $(i,j)$ to the four-element set of actions $A = \{ \leftarrow, \rightarrow, \uparrow, \downarrow \}$, i.e. an element of $A^X$.
    Initializing the value function $V$ with the environments immediate reward whose only non-zero value is $V(0,0)=1$ (top-left corner) and the policy with a upwards facing constant action $\pi(i,j)=\ \uparrow$ for all $(i,j)\in X$, a value improvement step would leave the policy unchanged while updating $V$ to $\overline{\pi}\comp\lambda\comp V$, which differs with $V$ only at $(0,1)\mapsto \beta$.

    If we instead perform a policy improvement step, the value function remains unchanged while the new policy differs with $\pi$ at $(1,0)\mapsto \arg\max_{a\in A}(\lambda\comp v)(1,0,a)=\ \leftarrow$.
\end{example}

\subsection{Continuous-space deterministic decision processes}

\begin{example}[Inverted pendulum]
    A state of our time-discretized inverted pendulum on a cart consists of $[y,\dot{y},\theta,\dot{\theta}]^\top$ in the state space $X=\R^4$.
    % Either the non-linear discrete-time setting, the propagator $F_{\Delta t}: x_k \to x_{k+1}$ in the non-linear setting,
    The linearized transition function that sends $x_k$ to $x_{k+1}=Ax_k+Ba_k$ is a smooth map $X\to Y$.
    % The integral expression of the reward or cost function is a smooth map that can be divided into a sum of integrals $\int_{t_k}^{t_{k+1}} C(x(t),u(t))dt$, but these individual integrals cannot be used as the backwards map because the input to the forwards map is a pair $x(t),u(t)$ for a specific $t$ value, while the required backwards cost is an integral of $C(x(t),u(t))$ between $[x(t_k),x(t_{k+1})]$. SAME argument for the forwards part actually.
    The discretized cost $J(x,a)=\sum_{k=0}^\infty \beta^k C(x(k),a(k))$ defines the backwards smooth function $X\times A\times\R\to \R$ which adds the cost at the $k$th time step $C(x(k),a(k))$ to the discounted sum:
   \[ \left(x(k),a(k),\sum_{j=k+1} \beta^jC(x(j),a(j))\right) \mapsto \sum_{j=k}\beta^j C(x(j),a(j)) \]
    These two maps form an optic $\lambda:\binom{X\times A}{\R}\to\binom{X}{\R}$ in $\Optic_\Euc$. Note that the cost function $C$ is itself typically not affine, but rather convex (intuitively, since the `good states' that should minimise the cost fall in the middle of the state space).

    In conclusion, the two optics involved in this example are
    \begin{equation*}
        \begin{aligned}[c]
            \lambda=\binom{f}{J}:\binom{X\times A}{\R}&\to\binom{X}{\R} \\
            f: X\times A &\to X \\
            (x,a)&\mapsto Ax + Ba \\
            J: X\times A \times \R &\to \R \\
            (x, a, r) &\mapsto C (x, a) + \beta r
        \end{aligned}
        \begin{aligned}[c]
            \overline\pi=\binom{\pi}{p_2}:\binom{X}{\R}&\to\binom{X\times A}{\R} \\
            \mathrm{gr}(\pi):X &\to X\times A \\
            x &\mapsto (x,\pi(x)) \\
            p_2:X\times \R&\to \R \\
            (\_,r) &\mapsto r
        \end{aligned}
    \end{equation*}
\end{example}

This formalisation of the continuous state space misses however a practical problem.
Let $S$ be a continuous state space.
In numerical implementations, policy improvement over $S$ needs to map an action to every point in the space.
Two common approaches are to discretize the state space into a possibly non-uniform grid, or to restrict the space of values to a family of parametrized functions \cite[Sec.4.]{rust_numerical_dp_in_econ}.
The discretization approach treats the continuous state space as a distribution over a simplicial complex $X$ obtained e.g. by triangulation, $\Euc(1,S)\approx \Mark(I,X)$, where a continuous state gets mapped to a distribution over the barycentric coordinates of the simplex.
This effectively transforms the initial continuous-space deterministic decision process into a discrete-space MDP, modelling numerical approximation errors as stochastic uncertainty.

\subsection{Discrete-space Markov decision processes}

Consider a Markov decision process with a discrete set of states $X$, discrete and finite set of actions $A$, a transition Markov kernel $f : X \times A \to \Delta (X)$, expected payoff function $U : X \times A \to \R$ and discount factor $\beta \in (0, 1)$. We can write the transition function as conditional probabilities $f (x' \mid x, a)$.

We can convert this data into an optic $\lambda : \binom{X \otimes A}{\R} \to \binom{X}{\R}$ in the category $\Optic_{\Mark, \Conv}$ given by $\Mark$ acting on both itself and $\Conv$. This optic is given concretely by $(X \otimes A, g, g')$ where $g : X \otimes A \to X \otimes A \otimes X$ in $\Mark$ is given by $\Delta_{X\otimes A}\comp(f\otimes \id_{X\otimes A})$, and $g' : \Delta (X \otimes A) \to [\R, \R]$ in $\Conv$ is defined by $g' (\alpha) (r) = \mathbb E U (\alpha) + \beta r$, where $\alpha \in \Delta (X \times A)$ is a joint distribution on states and actions. Alternatively, we can note that the domain of $g'$ is free on the set $X \times A$ (although it cannot be considered free on an object of $\Mark$), and define it as the linear extension of $g' (x, a) (r) = U (x, a) + \beta r$.

With this setup, value improvement $(\pi, V) \mapsto (\pi, \overline\pi\comp \lambda\comp V)$ yields the value function
\[ V' (x) = \E_{a \sim \pi (x)} [U (x, a) + \beta V (f (x, a))] \]
Alternating steps of value and policy improvement converge to the optimal policy $\pi^*$ and value function $V^*$, which maximises the expected value of the policy:
\[ V^*_{\pi^*} (x_0) = \E \sum_{k = 0}^\infty \beta^k U (x_k, \pi^* (x_k)) \]

\begin{example}[Gridworld, continued]
    In a proper MDP, transition functions can be stochastic, and update steps have to take expectations over values: value improvement maps $(\pi, V) \mapsto (\pi, \overline\pi\comp \lambda\comp V)$ and policy improvement maps $(\pi, V) \mapsto (x \mapsto \arg\max_{a \in A} \E (\lambda\comp V) (x, a), V)$.
%    Stochastic policy improvement steps like $\epsilon$-greedy can also be introduced in this framework.
    This model also accepts stochastic policy improvement steps like $\epsilon$-greedy, which is an ad hoc heuristic technique of balancing exploration and exploitation in reinforcement learning \cite[Sec.2]{kaelbling_rl_survey}, a problem which is known in control theory as the identification-control conflict.
\end{example}

\subsection{Continuous-space Markov decision processes}

For continuous-space MDPs we need a category of continuous Markov kernels. There are several possibilities for this arising as the kleisli category of a monad, such as the Giry monad on measurable spaces \cite{giry82}, the Radon monad on compact Hausdorff spaces \cite{swirszcz_monadic_functors_convexity} and the Kantorovich monad on complete metric spaces \cite{fritz_perrone_probability_monad_colimit}. However, control theorists typically work with more specific parametrised families of distributions for computational reasons, the most common being normal distributions. We will work with the category $\Gauss$ of Euclidean spaces and affine functions with Gaussian noise \cite[section 6]{Synthetic_approach}. (This is an example of a Markov category that does not arise as the kleisli category of a monad, because its multiplication map would not be affine.) This works because the pushforward measure of a Gaussian distribution along an affine function is still Gaussian, which fails for more general functions.

\begin{example}
	We will formulate the savings problem with normally-distributed income. The inequality constraints (namely that the balance cannot be negative and that the agent cannot consume more than their current balance) introduce nonlinearities. We can deal with the latter by constraining the optimisation in the policy improvement step, but the former threatens to take us outside the category $\Gauss$ and we must allow the balance to possibly become negative for the purposes of this example.
	
	$\Gauss$ is a Markov category that is not cartesian (the monoidal product is the cartesian product of Euclidean spaces, which adds the dimensions), so it acts on itself by the monoidal product and we take the category $\Optic_\Gauss$. We take the state and action spaces to be $X = A = \R$. The transition function $f : \R \otimes \R \to \R$ is given by $f (x, a) = (1 + \gamma) x - a + \mathcal N (\mu, \sigma)$, and the payoff function $U : \R \otimes \R \to \R$ is given by $U (x, a) = a$.
	
	We modify the policy improvement step to be
	\begin{align*}
		\text{\textbf{Policy improvement:}}& &&(\pi, V) \mapsto (x \mapsto \arg\max_{a \in A (x)} (\lambda\comp V) (x, a), V)
	\end{align*}
	where $A (x)$ is the set $A (x) = \{ a \in \R \mid 0 \leq a \leq x + i \}$. This enforces that the agent cannot consume negative amounts or consume more than their current balance - since the optimisation is done externally to the category $\Gauss$ we can avoid one source of nonlinearity this way.
	
	%For purposes of example, it is convenient to ignore the restriction that the bank balance cannot become negative, since the binary $\max$ operator is nonlinear and hence does not preserve Gaussians under pushforward measure. We therefore work in the category $\Gauss$, and take the state and action spaces to both be $X = A = \R$. (This technically causes the optimisation problem to become unbounded and degenerate, since there is no reason for the agent to not consume an arbitrarily high amount in every stage, but we ignore this problem.)

	%Since we do not have a corresponding Eilenberg-Moore category, the expectation over payoffs is not taken `automatically': a value improvement step $(\pi, V) \mapsto \overline\pi\comp \lambda\comp V$ will produce a function $V' : \R \to \R$ in $\Gauss$ with nonzero variance, even if $V$ always has zero variance. 
	%We can take an expectation by `manually' setting the variance to zero after every value improvement step. 
	%Finding a more appropriate categorical setting for this example (also allowing for the nonlinear $\min$ in payoffs to bound the optimisation problem) constitutes future work.
\end{example}

\section{Q-learning}
% \begin{align}
%     q_\pi(s,a) &= U(s,a) + \beta V_\pi(f(s,a)) \label{eq:q} \\
%     V_\pi'(s) &= q_\pi(s,\pi(s)) \\
%     \pi'(s) &= \arg\max_a q_\pi(s,a)
% \end{align}

Consider a deterministic decision process corresponding to the optic $\lambda : \binom{X \times A}{\R} \to \binom{X}{\R}$.
The dynamical system with state space $A^X\times \Optic_\Set\left(\binom{X\times A}{\R},I\right)$ has elements $(\pi,q)$ consisting of a \emph{state-action} value function $q: X\times A\to \R$  as in \eqref{eq:q} rather than a state-value function $V : X \to \R$.

We can define similar update steps
\begin{align*}
	\text{\textbf{Value improvement:}}& &&(\pi, q) \mapsto (\pi, \lambda\comp\overline\pi\comp q) \\
	\text{\textbf{Policy improvement:}}& &&(\pi, q) \mapsto \left( x \mapsto \arg\max_{a \in A} q (x, a), q \right)
\end{align*}
These can also be fused into a single step:
\begin{align*}
	\text{\textbf{State-action value iteration:}} &&&(\pi, q) \mapsto \left( x \mapsto \arg\max_{a \in A} q (x, a), \lambda\comp \overline\pi\comp q \right)
\end{align*}
Observe that composition of the $\lambda$ optic with $\overline\pi$ is flipped compared to the case seen in Section \ref{sec:4.1}, as we want an element of $\Optic_\Set\left(\binom{X\times A}{\R},\binom{X\times A}{\R}\right)$ to compose with $q$.

The advantage of learning state-action value functions $X \otimes A \to \R$ rather than state-value functions $X \to \R$ is that is gives a way to approximate $\arg\max_{a \in A} (\lambda\comp \overline\pi\comp q) (x, a)$ without making any use of $\lambda$, namely by instead using $\arg\max_{a \in A} q (x, a)$. This leads to an effective method known as Q-learning for computing optimal control policies even when the MDP is unknown, with only a single state transition and payoff being sampled at each time-step. This is the essential difference between classical control theory and \emph{reinforcement learning}. The above method, despite learning a $q$-function, is \emph{not} Q-learning because it makes use of $\lambda$ during value improvement.

Q-learning \cite{q_learning} is a sampling algorithm that approximates the state-action value iteration, usually by a lookup table $Q$ referred as Q-matrix.
It treats the optic as a black box, having therefore no access to the transition or rewards functions used in \eqref{eq:q}, and instead updates $q$ by interacting with the environment dynamics:
\[ q'(x',a) = (1-\alpha)q(x,a)+\alpha(r+\beta \max_{a'}q(x',a')) \]
where $\alpha\in(0,1)$ is a weighting parameter.
Note that both the new state $x'$ and the reward $r$ are obtained by interacting with the system, rather than looked ahead by $x'=f(x,a)$ and $r=U(x,a)$.
It falls in the family of temporal difference algorithms.

\section{Further work}

At the end of the previous section, it can be seen that Q-learning is no longer essentially using the structure of the category of optics, instead treating the Q-function as a mere function. We believe this can be overcome using the framework of categorical cybernetics \cite{towards_foundations_categorical_cybernetics}, leading to a fully optic-based approach to reinforcement learning. By combining with other instantiations of the same framework, it is hoped to encompass the zoo of modern variants of reinforcement learning that have achieved spectacular success in many applications in the last few years. For example, deep Q-learning represents the Q-function not as a matrix but as a deep neural network, trained by gradient descent, allowing much higher dimensionality to be handled in practice. Deep learning is currently one of the main applications of categorical cybernetics \cite{cruttwell_etal_categorical_foundations_gradient_based_learning}.

The proof that dynamic programming algorithms converge to the optimal policy and value function typically proceed by noting that the set of all value functions form a complete ordered metric space and that value improvement is a monotone contraction mapping. The metric structure is used to prove that iteration converges to a unique fixpoint by the contraction mapping theorem, and then the order structure is used to prove that this fixpoint is indeed optimal. Since value improvement is optic composition, these facts would be a special case of the category of optics being enriched in the category of ordered metric spaces and monotone contraction mappings. We do not currently know whether such an enrichment is possible. Unlike costates, general optics have nontrivial forwards passes, so there are two possible approaches: either ignore the forwards passes and defining a metric only in terms of the backwards passes, or defining a metric also using the forwards passes, for example using the Kantorovich metric between distributions. This would also be a reasonable place to unify our approach with the coalgebraic approach with metric coinduction \cite{LTVMDPs}.

Finally, continuous time MDPs pose a serious challenge to any approach for which categorical composition is sequencing in time, since composition of two morphisms in a category appears to be inherently discrete-time. (Open games are similarly unable to handle dynamic games with continuous time, such as pursuit games.) A plausible approach to this is to associate an endomorphism in a category to every real interval, by treating that interval of time as a single discrete time-step, and then requiring that all morphisms compose together correctly, similar to a sheaf condition. It is hoped that the Bellman-Jacobi-Hamilton equation, a PDE that is the continuous time analogue of the discrete-time Bellman equation, will similarly arise as a fixpoint in this way. Exploring this systematically is important future work.

%\clearpage
\bibliographystyle{eptcs}
\bibliography{Value_iteration}

\begin{thebibliography}{10}
\providecommand{\bibitemdeclare}[2]{}
\providecommand{\surnamestart}{}
\providecommand{\surnameend}{}
\providecommand{\urlprefix}{Available at }
\providecommand{\url}[1]{\texttt{#1}}
\providecommand{\href}[2]{\texttt{#2}}
\providecommand{\urlalt}[2]{\href{#1}{#2}}
\providecommand{\doi}[1]{doi:\urlalt{http://dx.doi.org/#1}{#1}}
\providecommand{\eprint}[1]{arXiv:\urlalt{https://arxiv.org/abs/#1}{#1}}
\providecommand{\bibinfo}[2]{#2}

\bibitemdeclare{article}{POMDP}
\bibitem{POMDP}
\bibinfo{author}{K.J \surnamestart {\AA}str{\"o}m\surnameend}
  (\bibinfo{year}{1965}): \emph{\bibinfo{title}{Optimal control of Markov
  processes with incomplete state information}}.
\newblock {\sl \bibinfo{journal}{Journal of Mathematical Analysis and
  Applications}} \bibinfo{volume}{10}(\bibinfo{number}{1}), pp.
  \bibinfo{pages}{174--205}, \doi{10.1016/0022-247x(65)90154-x}.

\bibitemdeclare{article}{baez_erbele_categories_control}
\bibitem{baez_erbele_categories_control}
\bibinfo{author}{John \surnamestart Baez\surnameend} \& \bibinfo{author}{Jason
  \surnamestart Erbele\surnameend} (\bibinfo{year}{2015}):
  \emph{\bibinfo{title}{Categories in control}}.
\newblock {\sl \bibinfo{journal}{Theory and applications of categories}}
  \bibinfo{volume}{30}(\bibinfo{number}{24}), pp. \bibinfo{pages}{836--881}.
\newblock \eprint{1405.6881}.

\bibitemdeclare{misc}{bakirtzis_categorical_semantics_rl}
\bibitem{bakirtzis_categorical_semantics_rl}
\bibinfo{author}{Georgios \surnamestart Bakirtzis\surnameend},
  \bibinfo{author}{Michail \surnamestart Savvas\surnameend} \&
  \bibinfo{author}{Ufuk \surnamestart Topcu\surnameend} (\bibinfo{year}{2022}):
  \emph{\bibinfo{title}{Categorical semantics of compositional reinforcement
  learning}}.
\newblock \eprint{2208.13687}.

\bibitemdeclare{phdthesis}{barfuss_learning_dynamics}
\bibitem{barfuss_learning_dynamics}
\bibinfo{author}{Wolfram \surnamestart Barfuss\surnameend}
  (\bibinfo{year}{2019}): \emph{\bibinfo{title}{Learning dynamics and decision
  paradigms in social-ecological dilemmas}}.
\newblock Ph.D. thesis, \bibinfo{school}{Humboldt-Universit\"at zu Berlin},
  \doi{10.18452/20127}.

\bibitemdeclare{article}{bellman_1957}
\bibitem{bellman_1957}
\bibinfo{author}{Richard \surnamestart Bellman\surnameend}
  (\bibinfo{year}{1957}): \emph{\bibinfo{title}{Dynamic Programming}}.
\newblock {\sl \bibinfo{journal}{Princeton University Press}},
  \doi{10.1515/9781400835386}.

\bibitemdeclare{book}{bertsekas_dpoc_vol2}
\bibitem{bertsekas_dpoc_vol2}
\bibinfo{author}{Dimitri~P. \surnamestart Bertsekas\surnameend}
  (\bibinfo{year}{2007}): \emph{\bibinfo{title}{Dynamic Programming and Optimal
  Control, Vol. II}}, \bibinfo{edition}{3rd} edition.
\newblock \bibinfo{publisher}{Athena Scientific}.

\bibitemdeclare{misc}{bolt_hedges_zahn_bayesian_open_games}
\bibitem{bolt_hedges_zahn_bayesian_open_games}
\bibinfo{author}{Joe \surnamestart Bolt\surnameend}, \bibinfo{author}{Jules
  \surnamestart Hedges\surnameend} \& \bibinfo{author}{Philipp \surnamestart
  Zahn\surnameend} (\bibinfo{year}{2019}): \emph{\bibinfo{title}{Bayesian open
  games}}.
\newblock \eprint{1910.03656}.
\newblock \bibinfo{note}{Forthcoming in \emph{Compositionality}}.

\bibitemdeclare{inproceedings}{botta13}
\bibitem{botta13}
\bibinfo{author}{Nicola \surnamestart Botta\surnameend}, \bibinfo{author}{Cezar
  \surnamestart Ionescu\surnameend} \& \bibinfo{author}{Edwin~C. \surnamestart
  Brady\surnameend} (\bibinfo{year}{2013}): \emph{\bibinfo{title}{Sequential
  decision problems, dependently-typed solutions}}.
\newblock In: {\sl \bibinfo{booktitle}{Joint Proceedings of the MathUI,
  OpenMath, {PLMMS} and ThEdu Workshops and Work in Progress at CICM, Bath,
  {UK}}}, {\sl \bibinfo{series}{{CEUR} Workshop Proceedings}}
  \bibinfo{volume}{1010}, \bibinfo{publisher}{CEUR-WS.org}.
\newblock \urlprefix\url{http://ceur-ws.org/Vol-1010/paper-06.pdf}.

\bibitemdeclare{misc}{braithwaite_etal_fibre_optics}
\bibitem{braithwaite_etal_fibre_optics}
\bibinfo{author}{Dylan \surnamestart Braithwaite\surnameend},
  \bibinfo{author}{Matteo \surnamestart Capucci\surnameend},
  \bibinfo{author}{Bruno \surnamestart Gavranovi\'c\surnameend},
  \bibinfo{author}{Jules \surnamestart Hedges\surnameend} \&
  \bibinfo{author}{Eigil~Fjeldgren \surnamestart Rischel\surnameend}
  (\bibinfo{year}{2022}): \emph{\bibinfo{title}{Fibre optics}}.
\newblock \eprint{2112.11145}.

\bibitemdeclare{misc}{capucci_gavranovic_actegories}
\bibitem{capucci_gavranovic_actegories}
\bibinfo{author}{Matteo \surnamestart Capucci\surnameend} \&
  \bibinfo{author}{Bruno \surnamestart Gavranovi\'c\surnameend}
  (\bibinfo{year}{2022}): \emph{\bibinfo{title}{Actegories for the working
  amthematician}}.
\newblock \eprint{2203.16351}.
\newblock \bibinfo{note}{ArXiv:2203.16351}.

\bibitemdeclare{article}{towards_foundations_categorical_cybernetics}
\bibitem{towards_foundations_categorical_cybernetics}
\bibinfo{author}{Matteo \surnamestart Capucci\surnameend},
  \bibinfo{author}{Bruno \surnamestart Gavranovi{\'c}\surnameend},
  \bibinfo{author}{Jules \surnamestart Hedges\surnameend} \&
  \bibinfo{author}{Eigil~Fjeldgren \surnamestart Rischel\surnameend}
  (\bibinfo{year}{2022}): \emph{\bibinfo{title}{Towards foundations of
  categorical cybernetics}}.
\newblock {\sl \bibinfo{journal}{Electronic Proceedings in Theoretical Computer
  Science}} \bibinfo{volume}{372}, pp. \bibinfo{pages}{235--248},
  \doi{10.4204/eptcs.372.17}.

\bibitemdeclare{article}{capucci_etal_translating_extensive_form}
\bibitem{capucci_etal_translating_extensive_form}
\bibinfo{author}{Matteo \surnamestart Capucci\surnameend},
  \bibinfo{author}{Neil \surnamestart Ghani\surnameend},
  \bibinfo{author}{J{\'{e}}r{\'{e}}my \surnamestart Ledent\surnameend} \&
  \bibinfo{author}{Fredrik~Nordvall \surnamestart Forsberg\surnameend}
  (\bibinfo{year}{2022}): \emph{\bibinfo{title}{Translating extensive form
  games to open games with agency}}.
\newblock {\sl \bibinfo{journal}{Electronic Proceedings in Theoretical Computer
  Science}} \bibinfo{volume}{372}, pp. \bibinfo{pages}{221--234},
  \doi{10.4204/eptcs.372.16}.

\bibitemdeclare{misc}{roman_etal_profunctor_optics}
\bibitem{roman_etal_profunctor_optics}
\bibinfo{author}{Bryce \surnamestart Clarke\surnameend}, \bibinfo{author}{Derek
  \surnamestart Elkins\surnameend}, \bibinfo{author}{Jeremy \surnamestart
  Gibbons\surnameend}, \bibinfo{author}{Fosco \surnamestart
  Loregian\surnameend}, \bibinfo{author}{Bartosz \surnamestart
  Milewski\surnameend}, \bibinfo{author}{Emily \surnamestart
  Pillmore\surnameend} \& \bibinfo{author}{Mario \surnamestart
  Rom\'an\surnameend} (\bibinfo{year}{2020}): \emph{\bibinfo{title}{Profunctor
  optics: {A} categorical update}}.
\newblock \eprint{2001.07488}.

\bibitemdeclare{misc}{cruttwell_etal_categorical_foundations_gradient_based_learning}
\bibitem{cruttwell_etal_categorical_foundations_gradient_based_learning}
\bibinfo{author}{G.S.H. \surnamestart Cruttwell\surnameend},
  \bibinfo{author}{Bruno \surnamestart Gavranovi\'c\surnameend},
  \bibinfo{author}{Neil \surnamestart Ghani\surnameend}, \bibinfo{author}{Paul
  \surnamestart Wilson\surnameend} \& \bibinfo{author}{Fabio \surnamestart
  Zanasi\surnameend} (\bibinfo{year}{2021}): \emph{\bibinfo{title}{Categorical
  foundations of gradient-based learning}}.
\newblock \eprint{2103.01931}.

\bibitemdeclare{article}{denardo_contraction_in_dp}
\bibitem{denardo_contraction_in_dp}
\bibinfo{author}{Eric~V. \surnamestart Denardo\surnameend}
  (\bibinfo{year}{1967}): \emph{\bibinfo{title}{Contraction Mappings in the
  Theory Underlying Dynamic Programming}}.
\newblock {\sl \bibinfo{journal}{{SIAM} Review}}
  \bibinfo{volume}{9}(\bibinfo{number}{2}), pp. \bibinfo{pages}{165--177},
  \doi{10.1137/1009030}.
\newblock
  \urlprefix\url{https://ui.adsabs.harvard.edu/abs/1967SIAMR...9..165D}.

\bibitemdeclare{misc}{eschenbaum_etal_robust_algorithmic_collusion}
\bibitem{eschenbaum_etal_robust_algorithmic_collusion}
\bibinfo{author}{Nicolas \surnamestart Eschenbaum\surnameend},
  \bibinfo{author}{Filip \surnamestart Mellgren\surnameend} \&
  \bibinfo{author}{Philipp \surnamestart Zahn\surnameend}
  (\bibinfo{year}{2022}): \emph{\bibinfo{title}{Robust algorithmic collusion}}.
\newblock \eprint{2201.00345}.

\bibitemdeclare{incollection}{LTVMDPs}
\bibitem{LTVMDPs}
\bibinfo{author}{Frank M.~V. \surnamestart Feys\surnameend},
  \bibinfo{author}{Helle~Hvid \surnamestart Hansen\surnameend} \&
  \bibinfo{author}{Lawrence~S. \surnamestart Moss\surnameend}
  (\bibinfo{year}{2018}): \emph{\bibinfo{title}{Long-{Term} {Values} in
  {Markov} {Decision} {Processes}, ({Co}){Algebraically}}}.
\newblock In \bibinfo{editor}{Corina \surnamestart C{\^\i}rstea\surnameend},
  editor: {\sl \bibinfo{booktitle}{Coalgebraic {Methods} in {Computer}
  {Science}}}, \bibinfo{volume}{11202}, \bibinfo{publisher}{Springer
  International Publishing}, \bibinfo{address}{Cham}, pp.
  \bibinfo{pages}{78--99}, \doi{10.1007/978-3-030-00389-0{\_}6}.
\newblock \urlprefix\url{http://link.springer.com/10.1007/978-3-030-00389-0_6}.
\newblock \bibinfo{note}{Series Title: Lecture Notes in Computer Science}.

\bibitemdeclare{book}{control_systems_design}
\bibitem{control_systems_design}
\bibinfo{author}{Bernard \surnamestart Friedland\surnameend}
  (\bibinfo{year}{1985}): \emph{\bibinfo{title}{Control Systems Design}}.
\newblock \bibinfo{publisher}{McGraw-Hill Companies}.

\bibitemdeclare{misc}{fritz09}
\bibitem{fritz09}
\bibinfo{author}{Tobias \surnamestart Fritz\surnameend} (\bibinfo{year}{2009}):
  \emph{\bibinfo{title}{Convex spaces {I}: {D}efinitions and examples}}.
\newblock \eprint{0903.5522}.

\bibitemdeclare{article}{Synthetic_approach}
\bibitem{Synthetic_approach}
\bibinfo{author}{Tobias \surnamestart Fritz\surnameend} (\bibinfo{year}{2020}):
  \emph{\bibinfo{title}{A synthetic approach to {Markov} kernels, conditional
  independence and theorems on sufficient statistics}}.
\newblock {\sl \bibinfo{journal}{Advances in Mathematics}}
  \bibinfo{volume}{370}, p. \bibinfo{pages}{107239},
  \doi{10.1016/j.aim.2020.107239}.
\newblock \eprint{1908.07021}.

\bibitemdeclare{article}{fritz_perrone_probability_monad_colimit}
\bibitem{fritz_perrone_probability_monad_colimit}
\bibinfo{author}{Tobias \surnamestart Fritz\surnameend} \&
  \bibinfo{author}{Paolo \surnamestart Perrone\surnameend}
  (\bibinfo{year}{2019}): \emph{\bibinfo{title}{A probability monad as the
  colimit of spaces of finite samples}}.
\newblock {\sl \bibinfo{journal}{Theory and applications of categories}}
  \bibinfo{volume}{34}(\bibinfo{number}{7}), pp. \bibinfo{pages}{170--220}.
\newblock \eprint{1712.05363}.

\bibitemdeclare{inproceedings}{hedges_etal_compositional_game_theory}
\bibitem{hedges_etal_compositional_game_theory}
\bibinfo{author}{Neil \surnamestart Ghani\surnameend}, \bibinfo{author}{Jules
  \surnamestart Hedges\surnameend}, \bibinfo{author}{Viktor \surnamestart
  Winschel\surnameend} \& \bibinfo{author}{Philipp \surnamestart
  Zahn\surnameend} (\bibinfo{year}{2018}): \emph{\bibinfo{title}{Compositional
  game theory}}.
\newblock In: {\sl \bibinfo{booktitle}{Proceedings of Logic in Computer Science
  (LiCS) 2018}}, \bibinfo{publisher}{{ACM}}, pp. \bibinfo{pages}{472--481},
  \doi{10.1145/3209108.3209165}.

\bibitemdeclare{incollection}{giry82}
\bibitem{giry82}
\bibinfo{author}{Mich\`elle \surnamestart Giry\surnameend}
  (\bibinfo{year}{1982}): \emph{\bibinfo{title}{A categorical approach to
  probability theory}}.
\newblock In: {\sl \bibinfo{booktitle}{Lecture Notes in Mathematics}},
  \bibinfo{publisher}{Springer Berlin Heidelberg}, pp. \bibinfo{pages}{68--85},
  \doi{10.1007/bfb0092872}.

\bibitemdeclare{misc}{hedges_coherence_lenses_open_games}
\bibitem{hedges_coherence_lenses_open_games}
\bibinfo{author}{Jules \surnamestart Hedges\surnameend} (\bibinfo{year}{2017}):
  \emph{\bibinfo{title}{Coherence for lenses and open games}}.
\newblock \eprint{1704.02230}.

\bibitemdeclare{inproceedings}{hedges_morphisms_open_games}
\bibitem{hedges_morphisms_open_games}
\bibinfo{author}{Jules \surnamestart Hedges\surnameend} (\bibinfo{year}{2018}):
  \emph{\bibinfo{title}{Morphisms of open games}}.
\newblock In: {\sl \bibinfo{booktitle}{Proceedings of {MFPS} 2018}},
  \bibinfo{volume}{341}, \bibinfo{publisher}{Elsevier {BV}}, pp.
  \bibinfo{pages}{151--177}, \doi{10.1016/j.entcs.2018.11.008}.

\bibitemdeclare{article}{kaelbling_rl_survey}
\bibitem{kaelbling_rl_survey}
\bibinfo{author}{L.~P. \surnamestart Kaelbling\surnameend},
  \bibinfo{author}{M.~L. \surnamestart Littman\surnameend} \&
  \bibinfo{author}{A.~W. \surnamestart Moore\surnameend}
  (\bibinfo{year}{1996}): \emph{\bibinfo{title}{Reinforcement Learning: A
  Survey}}.
\newblock {\sl \bibinfo{journal}{Journal of Artificial Intelligence Research}}
  \bibinfo{volume}{4}, pp. \bibinfo{pages}{237--285}, \doi{10.1613/jair.301}.

\bibitemdeclare{article}{kock_closed_categories_commutative_monads}
\bibitem{kock_closed_categories_commutative_monads}
\bibinfo{author}{Anders \surnamestart Kock\surnameend} (\bibinfo{year}{1971}):
  \emph{\bibinfo{title}{Closed categories generated by commutative monads}}.
\newblock {\sl \bibinfo{journal}{Journal of the Australian Mathematical
  Society}} \bibinfo{volume}{12}(\bibinfo{number}{4}), pp.
  \bibinfo{pages}{405--424}, \doi{10.1017/s1446788700010272}.

\bibitemdeclare{book}{ljungqvist_sargent_recursive_macroeconomic_theory}
\bibitem{ljungqvist_sargent_recursive_macroeconomic_theory}
\bibinfo{author}{Lars \surnamestart Ljungqvist\surnameend} \&
  \bibinfo{author}{Thomas \surnamestart Sargent\surnameend}
  (\bibinfo{year}{2004}): \emph{\bibinfo{title}{Recursive macroeconomic
  theory}}.
\newblock \bibinfo{publisher}{MIT Press}.

\bibitemdeclare{book}{loregian-coend-cofriend}
\bibitem{loregian-coend-cofriend}
\bibinfo{author}{Fosco \surnamestart Loregian\surnameend}
  (\bibinfo{year}{2021}): \emph{\bibinfo{title}{Coend calculus}}.
\newblock \bibinfo{publisher}{Cambridge University Press}.
\newblock \eprint{1501.02503}.

\bibitemdeclare{article}{mody_steffen_optimal_charging_electric_vehicles}
\bibitem{mody_steffen_optimal_charging_electric_vehicles}
\bibinfo{author}{Sagar \surnamestart Mody\surnameend} \&
  \bibinfo{author}{Thomas \surnamestart Steffen\surnameend}
  (\bibinfo{year}{2015}): \emph{\bibinfo{title}{Optimal charging of electric
  vehicles using a stochastic dynamic programming model and price prediction}}.
\newblock {\sl \bibinfo{journal}{{SAE} International Journal of Passenger Cars
  - Electronic and Electrical Systems}}
  \bibinfo{volume}{8}(\bibinfo{number}{2}), pp. \bibinfo{pages}{379--393},
  \doi{10.4271/2015-01-0302}.

\bibitemdeclare{article}{DJM}
\bibitem{DJM}
\bibinfo{author}{David~Jaz \surnamestart Myers\surnameend}
  (\bibinfo{year}{2021}): \emph{\bibinfo{title}{Double {Categories} of {Open}
  {Dynamical} {Systems} ({Extended} {Abstract})}}.
\newblock {\sl \bibinfo{journal}{Electronic Proceedings in Theoretical Computer
  Science}} \bibinfo{volume}{333}, pp. \bibinfo{pages}{154--167},
  \doi{10.4204/eptcs.333.11}.
\newblock \urlprefix\url{http://arxiv.org/abs/2005.05956v2}.

\bibitemdeclare{book}{Puterman}
\bibitem{Puterman}
\bibinfo{author}{Martin~L. \surnamestart Puterman\surnameend}
  (\bibinfo{year}{2005}): \emph{\bibinfo{title}{Markov decision processes:
  discrete stochastic dynamic programming}}.
\newblock \bibinfo{series}{Wiley series in probability and statistics},
  \bibinfo{publisher}{Wiley-Interscience}, \bibinfo{address}{Hoboken, NJ}.

\bibitemdeclare{misc}{Riley}
\bibitem{Riley}
\bibinfo{author}{Mitchell \surnamestart Riley\surnameend}
  (\bibinfo{year}{2018}): \emph{\bibinfo{title}{Categories of {Optics}}}.
\newblock \eprint{1809.00738}.

\bibitemdeclare{incollection}{rust_numerical_dp_in_econ}
\bibitem{rust_numerical_dp_in_econ}
\bibinfo{author}{John \surnamestart Rust\surnameend} (\bibinfo{year}{1996}):
  \emph{\bibinfo{title}{Numerical dynamic programming in economics}}.
\newblock In: {\sl \bibinfo{booktitle}{Handbook of Computational Economics}},
  \bibinfo{publisher}{Amsterdam: Elsevier}, pp. \bibinfo{pages}{619--729},
  \doi{10.1016/s1574-0021(96)01016-7}.

\bibitemdeclare{article}{Spivak}
\bibitem{Spivak}
\bibinfo{author}{David~I. \surnamestart Spivak\surnameend}
  (\bibinfo{year}{2020}): \emph{\bibinfo{title}{Generalized {Lens} {Categories}
  via functors {C}{\textasciicircum}op-{\textgreater}{Cat}}}.
\newblock \eprint{1908.02202}.

\bibitemdeclare{misc}{sturtz_categorical_probability}
\bibitem{sturtz_categorical_probability}
\bibinfo{author}{Kirk \surnamestart Stirtz\surnameend} (\bibinfo{year}{2015}):
  \emph{\bibinfo{title}{Categorical probability theory}}.
\newblock \eprint{1406.6030}.

\bibitemdeclare{book}{stokey_lucas_recursive_methods_economic_dynamics}
\bibitem{stokey_lucas_recursive_methods_economic_dynamics}
\bibinfo{author}{Nancy \surnamestart Stokey\surnameend},
  \bibinfo{author}{Robert \surnamestart Lucas\surnameend} \&
  \bibinfo{author}{Edward \surnamestart Prescott\surnameend}
  (\bibinfo{year}{1989}): \emph{\bibinfo{title}{Recursive methods in economic
  dynamics}}.
\newblock \bibinfo{publisher}{Harvard University Press},
  \doi{10.2307/j.ctvjnrt76}.

\bibitemdeclare{book}{RL_Intro}
\bibitem{RL_Intro}
\bibinfo{author}{Richard~S. \surnamestart Sutton\surnameend} \&
  \bibinfo{author}{Andrew~G. \surnamestart Barto\surnameend}
  (\bibinfo{year}{2020}): \emph{\bibinfo{title}{Reinforcement {Learning}: {An}
  {Introduction}}}.
\newblock \bibinfo{publisher}{MIT Press}.

\bibitemdeclare{article}{swirszcz_monadic_functors_convexity}
\bibitem{swirszcz_monadic_functors_convexity}
\bibinfo{author}{T.~\surnamestart \'{S}wirszcz\surnameend}
  (\bibinfo{year}{1974}): \emph{\bibinfo{title}{Monadic functors and
  convexity}}.
\newblock {\sl \bibinfo{journal}{Bulletin de l'Aca\'{e}demie Polonaise des
  Sciences}} \bibinfo{volume}{22}(\bibinfo{number}{1}).

\bibitemdeclare{misc}{vertechi_dependent_optics}
\bibitem{vertechi_dependent_optics}
\bibinfo{author}{Pietro \surnamestart Vertechi\surnameend}
  (\bibinfo{year}{2022}): \emph{\bibinfo{title}{Dependent optics}}.
\newblock \eprint{2204.09547}.

\bibitemdeclare{article}{q_learning}
\bibitem{q_learning}
\bibinfo{author}{Christopher J. C.~H. \surnamestart Watkins\surnameend} \&
  \bibinfo{author}{Peter \surnamestart Dayan\surnameend}
  (\bibinfo{year}{1992}): \emph{\bibinfo{title}{Q-learning}}.
\newblock {\sl \bibinfo{journal}{Machine Learning}}
  \bibinfo{volume}{8}(\bibinfo{number}{3-4}), pp. \bibinfo{pages}{279--292},
  \doi{10.1007/bf00992698}.

\end{thebibliography}

\end{document}